\documentclass{article}
\usepackage{amsfonts,amsmath}
        \newtheorem{theorem}{Theorem}

        \newtheorem{remark}[theorem]{Remark}

\usepackage{epsfig}
\title{An iterative  domain decomposition method for free boundary problems with nonlinear flux jump
constraint}

\author{Juan Galvis\thanks{ISC/IAMCS, Texas A \& M University, College Station, TX 77843
({\tt jugal@math.tamu.edu})} \and H. M. Versieux \thanks{Instituto de Matem\'atica, Universidade Federal do Rio de Janeiro, Rio de Janeiro RJ ({\tt henrique@im.ufrj.br})}
}

\begin{document}
\maketitle



\section*{Abstract}
In this paper we design an iterative domain decomposition  method   for
free boundary problems with  nonlinear flux jump condition. Our approach is related to
 damped Newton's methods. The proposed scheme  requires, in each iteration,
the approximation of the flux on (both sides of) the free interface.
We present a Finite Element implementation of our method.
The numerical
implementation uses harmonically deformed triangulations
to inexpensively generate  finite element meshes in  subdomains.
We apply our method to a simplified model
for jet flows in pipes  and  to a simple magnetohydrodynamics model. Finally, we present numerical examples studying the  convergence of our scheme.
\section{Introduction}\label{sec:introduction}

In this paper, we propose a numerical iterative method for approximating  the solutions of  free boundary problems
in two dimensions. Our  iterative method for free boundary problems is based on
Domain Decomposition and damped Newton's method ideas.
In general terms, free boundary problems seek to determine unknown function $u$
with some prescribed conditions on a unknown interior interface, exterior boundary or (sub)domain.
In many applications, it is prescribed the value of $u$ on the free interface and it is required
that $u$ satisfy a condition involving (both  sides) derivatives of $u$ on the interface.
We mention jump conditions of Stefan, Bernoulli and Gibbs-Thomson type, among others.
There is a considerable literature of iterative methods for these type of free boundary problems; see for instance \cite{MR2149216,eppler,MR1450755,MR1677147,MR2306261,MR2639606,MR1427459} and references therein. In particular,  numerical finite elements methods have been proposed to solve Stefan-like
free boundary  problems
(including time dependent problems) and some other similar phase transition problems; see for instance
\cite{MR789678,Che-Shi-Yue,MR1079028,MR1114983,MR1094958}.  These methods use a variational formulation of their original problem.  Level set approach for Stefan problems
 were also proposed in \cite{MR1461705} and references therein.

The free boundary conditions that we deal with, up to our knowledge, have not
been extensively studied from the numerical point of view.
We are particularly interested in free boundary problems were
the unknown function satisfy nonlinear jump constraints across the free interface. More precisely,
given $\Omega\subset \mathbb{R}^2$, $g:\partial\Omega\to\mathbb{R}$ and
$\lambda\in \mathbb{R}$, we want to find a function
  $u:\Omega\to \mathbb{R}$ and a free interface $\Gamma$ (diving
$\Omega$ in two subdomains $\Omega^+=\{u>0\}$
and $\Omega^-=\{u<0\}$) such that $u$ and $\Gamma$ satisfy the following subdomain
equations and boundary condition,
\begin{eqnarray}\label{eq:problema-orig-ACF_general}
-\mbox{div}(a_1\nabla u)=0 & \mbox{ in }~\Omega^+,\\
-\mbox{div}(a_2\nabla u)=0 & \mbox{ in }~\Omega^-,\label{eq:problema-orig-ACF_general2}\\
    u=g& \mbox{ on}~ \partial \Omega \label{eq:problema-orig-ACF_general3}
\end{eqnarray}
and free interface condition
\begin{equation}\label{eq:problema-orig-ACF_generalFreeBoundary}
 a_1|\nabla u^+|^2  - a_2|\nabla u^-|^2= \lambda    ~\mbox{on} ~\Gamma,
\end{equation}
or, similar nonlinear constraint for the jump in the  derivative of $u$ across the free interface $\Gamma$. Above, $u^+$ and $u^-$ denote the value of
the solution $u$ on both sides of the free interface and the derivatives and the quantities
involved are interpreted as side limits. In many applications, the interface conditions
are imposed in a weak sense. These conditions can be also interpreted if we replace
the operators involved (e.g.,  trace of the  derivatives) by some smooth
or regularize version of them when necessary.

We are not aware of a simple inexpensive numerical  method  to solve problem
(\ref{eq:problema-orig-ACF_general})-(\ref{eq:problema-orig-ACF_generalFreeBoundary}).
The finite element  methods mentioned earlier to handle Stefan, Bernoulli and similar free boundary conditions are based on variational formulations. They do not seem to be easily extended to handle our nonlinear free boundary constraint.  Also,  Bernoulli type free boundary problems
 when one of the phases is a constant function seems to be easier to handle numerically. In this case, using
the fact that the tangential derivative on the free interface is zero and that the flux sign
can be a priori determined, the interface condition reduces to a linear condition
of the form $\partial_\eta u=\lambda$
where $\partial_\eta$ is the  normal derivative on the free interface.

We have two main applications in mind: 1)  the jet flow model studied by Alt, Caffarelli, Friedman \cite{Alt-Caffarelli-Friedman,Alt-Caf-Fri_Jet-I} and 2) a free boundary problem arising in magnetohydrodynamics
studied in \cite{friedman1,plasma2}. These applications are simplified mathematical versions of
complicated flow models and they focus in the main modeling aspects.
Despite of the mathematical simplifications, in either case, the resulting model problem above is still complex and finding and understanding solutions requires
numerical methods. The methods used for this problems  should be inexpensive and simple.
The method presented here is designed having these considerations into account. It can also be easily extended to handle different free boundary problems such as the stationary solutions of the Stefan's problem, and other similar problems.

The iterative method proposed in this paper for problem
(\ref{eq:problema-orig-ACF_general})-(\ref{eq:problema-orig-ACF_generalFreeBoundary}) is based on the following simple ideas. Assume the solution $u$ is sufficiently regular, and let $\Gamma$ denote the free boundary of problem  (\ref{eq:problema-orig-ACF_general}).  Since $u=0$ on $\Gamma$, $\nabla u^+= \partial_\eta u^+ \eta$ on $\Gamma$, where  $\eta$ is the outer normal vector of the region defined by the support of $u^+$. Hence, the free boundary condition  (\ref{eq:problema-orig-ACF_generalFreeBoundary}) reads
\begin{equation}\label{cond:front-livre}
 a_1|\partial_\eta u^+|^2  - a_2|\partial_\eta u^-|^2= \lambda    ~\mbox{on} ~\Gamma.
\end{equation}
Next, assume we have an approximation  $\widetilde{\Gamma}$ of $\Gamma$ dividing $\Omega$ in two different regions $\widetilde{\Omega}^+$ and $\widetilde{\Omega}^-$. We also assume that $\widetilde{\Omega}^+$ and $\widetilde{\Omega}^-$ are connected subdomains, and
$\partial \widetilde{\Omega}^+ \cap \partial \Omega=\Sigma^+$ and
 $\partial \widetilde{\Omega}^- \cap \partial \Omega=\Sigma^-$.
 In order to construct an approximation $\widetilde{u}$ of $u$,  we can  solve Dirichlet problems
(\ref{eq:problema-orig-ACF_general})-(\ref{eq:problema-orig-ACF_general3})
in the approximated subdomains with homogeneous Dirichlet boundary condition on
the approximated free interface $\widetilde{\Gamma}$.
The solution of these two independent problems give  $\widetilde{u}^+$
and $\widetilde{u}^-$.
   We observe that we do not expect the function $\widetilde{u}$ to satisfy condition (\ref{cond:front-livre}),
since $\widetilde{\Gamma}$ is only an approximation of $\Gamma$. Finally, we  update the approximation
of the free boundary by using the
quantity $\sigma=a_1|\partial_\eta \widetilde{u}^+(x)|^2  - a_2|\partial_\eta \widetilde{u}^-(x)|^2 - \lambda$
and a perturbation of $\widetilde{\Gamma}$ in its normal direction $\eta(x)$.
More specifically, we locally move $\widetilde{\Gamma}$ in the direction of $\eta(x)$
by a magnitude $\tau\sigma$ where $\tau$ is a positive damping parameter.
 Once the new approximation of $\Gamma$ is obtained we restart this procedure.

The rest of the paper is organized as follows. In Section \ref{sec:iterativo} we describe our iterative scheme. Section \ref{sec:fem} describes the finite element implementation of our method. In Section
\ref{sec:jet} we present the jet flow model proposed by  Alt, Caffarelli, Friedman and some numerical solutions for this problem. Numerical experiments for the  magnetohydrodynamics problem
studied in \cite{friedman1,plasma2} are presented in Section \ref{sec:mhd}. Section \ref{sec:additional-num}  presents some numerical experiments where we study convergence properties of our scheme. Finally, we present our conclusions and comments  in Section \ref{sec:conclusions}.

\section{Model problem and iterative method for the
free interface}\label{sec:iterativo}

In order to simplify  the presentation and fix ideas, we consider
the two dimensional case $\Omega \subset \mathbb{R}^2$, and the following free boundary problem
\begin{equation}\label{eq:differential-problem}
\left\{ \begin{array}{rll}
-\nabla \cdot(a_1\nabla u^{+})=&0& \mbox{ in } \{u>0\}\\
-\nabla \cdot(a_2\nabla u^{-})=&0& \mbox{ in } \{u<0\}\\
a_1(\partial_{\eta} u^{+} )^2-a_2(\partial_{\eta} u^{-})^2=&
\lambda  & \mbox{ on } \Gamma=\{u=0\}  \\
u=&g& \mbox{ on }  \partial \Omega.
\end{array}
\right.
\end{equation}
We also assume there exist two connected curves  $\Sigma^+, ~\Sigma^-,$ ~ such that
$\partial \Omega= \bar{\Sigma^+} \cup  \bar{\Sigma^-}$, and $g|_{\Sigma^+}>0$ and $g|_{\Sigma^-}<0$.
This model problem, or similar system of equations, appear in different applications.

We approximate the
solution  of problem
(\ref{eq:differential-problem}) by constructing
a sequence of   approximations of the free boundary $\Gamma$.  Assume we have an approximation of the free boundary $\Gamma$, then we solve two independent elliptic problems and
use the condition  $a_1(\partial_{\eta} u^{+} )^2-a_2(\partial_{\eta} u^{-})^2=\lambda$ to updated the
approximation of the free boundary as follows.

Assume $\Gamma_{n}$ is an approximation of the free boundary dividing the
domain $\Omega$ into two subdomains, $\Omega^+_n$ (enclosed by $\Sigma^+ \cup \Gamma_n$)
and $\Omega^-_n$ (enclosed by $\Sigma^- \cup \Gamma_n$).
We define the $n$-th approximation of $u$ as follows.   In $\Omega^+_n$
the function $u_n$ solves,
\begin{equation}\label{eq:current-u-plus}
\left\{ \begin{array}{rll}
-\nabla \cdot(a_1\nabla u_n)=&0& \mbox{ in } \Omega^+_n\\
u_n=&g& \mbox{ on }  \Sigma^+\\
u_n=&0& \mbox{ on }   \Gamma_n.\\
\end{array}
\right.
\end{equation}
In $\Omega^-_n$ the function $u_n$ solves,
\begin{equation}\label{eq:current-u-minus}
\left\{ \begin{array}{rll}
-\nabla \cdot(a_2\nabla u_n)=&0& \mbox{ in } \Omega^-_n\\
u_n=&g& \mbox{ on }  \Sigma^-\\
u_n=&0& \mbox{ on }  \Gamma_n.\\
\end{array}
\right.
\end{equation}

The main idea to define the updated approximation of the free boundary $\Gamma_{n+1}$  is very simple.
First, we define
\begin{equation}\label{eq:def:sigma}
\sigma = a_1(\partial_\eta u_n^{+} )^2-a_2( \partial_\eta u_n^{-})^2- \lambda.
\end{equation}
Here, we use the notation  $\partial_\eta u^\pm_n$ as the outward normal derivative of $u^\pm_n$  with respect to the region $\Omega^\pm_n$.
Next, if for instance, $\sigma(x)>0$ for some point $x\in \Gamma_n$,  then we would like  to locally update $\Gamma_n$ such that $\sigma(x)$ is closer to zero.
This can be done by   decreasing  the flux of $u_n$ in $\Omega^+_n$ and/or increasing  the flux of $u_n$ in $\Omega^-_n$
in a neighborhood of that point.  We expect to obtain this  by locally moving the free interface $\Gamma_n$ in the normal direction outward to $\Omega^+_n$.  We define the new approximation of the free interface
by
\begin{equation}\label{eq:next-Gamma}
\Gamma_{n+1}=\{x+\tau \sigma \vec{\eta}_{\Gamma^+_{n}} ~~; ~~
\mbox{ with } x\in \Gamma_{n}\}.
\end{equation}
Here $\tau=\tau(\sigma)$ is a small positive parameter, and $\vec{\eta}_{\Gamma^+_{n}}$ represents the unitary normal vector of $\Gamma_n$  outward to $\Omega^+_n$.

Finally, we observe that there are several ways to define $\Gamma_0$ dividing the domain $\Omega$ in two parts as desired. For instance, we can take $\Gamma_0$ as the
zero level set of any regular extension of the boundary  data $g$.

\begin{remark}
We note that we need only an approximation
of $~\sigma$ (which requires
only approximation of the flux). This is important in
case $u$ is not regular enough to
allow the computation of the square
of the flux.
\end{remark}

\begin{remark}
We mention that in \cite{Alt-Caffarelli-Friedman} it is proved that the solution of problem (\ref{eq:problema-orig-ACF_general})-(\ref{eq:problema-orig-ACF_generalFreeBoundary}), in
the case $a_1=a_2=1$, is a minimizer of the following functional
$ J(v)=\int_{\Omega} 1_{v>0}\left(
|\nabla v|^2+\lambda_1^{2}\right)+
1_{v<0}\left( |\nabla v|^2+\lambda_2^{2}\right)\;dx $
where $\lambda= \lambda_1^{2} - \lambda_2^{2}$, and $1_{v>0}$  is the characteristic function
of the set $[\{v>0\}:=\{x\in \Omega, v(x)>0\}]$, (similar for $1_{v<0}$). We have also developed  a method for this problem based on the minimization of this functional. This was performed by,  first, introducing a regularized approximation $J_\epsilon$ of the functional $J$. Next, we looked for  a minimum of  the functional $J_\epsilon$ by solving the steepest descent evolution PDE associated to  this functional.  However, the observed numerical results were not satisfactory. We also observe that this method requires to solve a nonlinear problem for each time step resulting in more computational work compared to our iterative method.
\end{remark}

We also observe that our  method  to solve problem (\ref{eq:differential-problem}) can also handle different problems. For instance, the same ideas apply to the following
abstract free boundary problem.
Let ${\cal L}_+$ and ${\cal L}_-$ represent two second order elliptic operators.
Assume $\Omega \subset \mathbb{R}^n$, and let $\phi$ and $\psi: \mathbb{R}\rightarrow \mathbb{R}$
 be two increasing functions.  Consider the problem of finding $u$ and a
free interface $\Gamma$ such that
\begin{equation}\label{eq:problema-geral}\left\{
 \begin{array}{lll}
 {\cal L}_+  u=0 & \mbox{ in }~\Omega^+,  \\
 {\cal L}_-  u=0 & \mbox{ in }~\Omega^-,  \\
  \phi(\partial_\eta u^+)   -  \psi(\partial_\eta u^-) = \lambda(x)  &  \mbox{on the free boundary}~ \Gamma\\
    u=g& \mbox{ on}~ \partial \Omega.
\end{array}\right.
\end{equation}
Here $\Gamma$ represents the free boundary separating the two phases,  $\partial_\eta u_i$ represents the outward normal derivative with respect to the i-th phase.
Up to our knowledge there is no
rigorous studies of such general class  of problems. We mention that these problems include
the stationary solutions of two phase Stefan problem  (see \cite{MR0141895,MR0227625,MR0351348,MR2508170} and references therein); and  other problems involving nonlinear free boundary conditions (see \cite{lei-que-teix}).

\section{Finite element implementation}\label{sec:fem}
Now we describe the finite element implementation
of our iterative method.
In each iteration we have to approximate the solutions
of  problems
(\ref{eq:current-u-plus}) and
(\ref{eq:current-u-minus}) as well as
$\sigma$ in (\ref{eq:def:sigma}).

Let $n \in \mathbb{N}$ be our iteration parameter,  $\mathcal{T}_n^h$ be a triangulation of $\Omega$ with nodes
$\{x_{n,j}^h\}_{j=1}^{N_{v}}$ and edges $\{e_{n,\ell}\}_{\ell=1}^{N_e}$, and let $\Gamma_n^h$  be an  approximation of the free boundary $\Gamma$, such that  $\Gamma_n^h \subset \cup_{\ell=1}^{N_e} e_{n,\ell}$.
Here we also assume $\Gamma_{n}^h$ divides  $\Omega$ into two subdomains, $\Omega^+_{n,h}$ (enclosed by $\Sigma^+ \cup \Gamma_n^h$)
and $\Omega^-_{n,h}$ (enclosed by $\Sigma^- \cup \Gamma_n^h$). Set $V^{+,h}_0= \{v\in {\cal P}^1(\mathcal{T}^h_n, \Omega^+_{n,h}); ~ v=0 ~\mbox{on}~\partial \Omega^+_{n,h} \}$, where
${\cal P}^1(\mathcal{T}^h_n, \Omega^+_{n,h})$ represents the set of continuous piecewise linear functions on $\mathcal{T}_n$ (the space $V^{-,h}_0$ is defined similarly). The n-th
 approximation of the solution $u$ of (\ref{eq:differential-problem}), denoted by
$u^h_n$, solves the finite element problems,
\begin{equation}\label{eq:current-u-plus-h}
\left\{ \begin{array}{rll}
\int_{\Omega^+_{n,h}} a_1\nabla u_n^h\nabla z^h\;dx=&0& \mbox{ for all }
z^h\in V^{+,h}_0 \\
u^h_n(x^h)=&g(x^h)& \mbox{ for all }  x^h\in \Sigma^+ , ~\mbox{and}~x^h\in \mathcal{T}^h_n\\
u^h(x^h)=&0& \mbox{ for all  }   x^h\in \Gamma_n^h , ~\mbox{and}~x^h\in \mathcal{T}^h_n,\\
\end{array}
\right.
\end{equation}
and
\begin{equation}\label{eq:current-u-minus-h}
\left\{ \begin{array}{rll}
\int_{\Omega^-_{n,h}} a_2\nabla u_n^h\nabla z^h\;dx =&0& \mbox{ for all }
z^h \in V^{-,h}_0\\
u^h_n(x^h)=&g(x^h)& \mbox{ for all }  x^h\in \Sigma^-, ~\mbox{and}~x^h\in \mathcal{T}^h_n\\
u^h(x^h)=&0& \mbox{ for all  }   x^h\in \Gamma_n^h  , ~\mbox{and}~x^h\in \mathcal{T}^h_n.\\
\end{array}
\right.
\end{equation}

We define  $\rho^+_{n,h}$ as an appropriate piecewise linear  approximation of the  flux
of $u_n^h|_{\Omega^+_{n,h}}$ across $\Gamma_n^h$;  see Appendix \ref{sec:flux}.
Analogously, we define   $\rho^-_{n,h}$ as the discrete flux
of $u_n^h|_{\Omega^-_{n,h}}$ across $\Gamma_n^h$.
We note that,
\[
\rho^+_{n,h}=\sum_{x_{n,j}^h\in \Gamma_n^h}
\alpha_{n,j}^+\psi_{j} ~~\mbox{and}~~\rho^-_{n,h}=\sum_{x_{n,j}^h\in \Gamma_n^h}
\alpha_{n,j}^-\psi_{j},
\]
where the function $\psi_j$ represents  a basis for the space ${\cal P}^1(\mathcal{T}_n, \Omega)$  restricted to $\Gamma^h$.

Define,
\begin{equation}\label{eq:def:sigmah}
\sigma^h_n = \sum_{x_{n,j}^h\in \Gamma_n^h}
\Big( (\alpha_{n,j}^+)^2-(\alpha_{n,j}^-)^2-\lambda^2 \Big)
\psi_{j}.
\end{equation}
The new approximation of the free interface is given by
the piecewise linear curve,
\begin{equation}\label{eq:next-Gamma-h}
\Gamma_{n+1}^h=\{x+\tau \sigma^h_n(x) \vec{\eta}_{\Gamma_{n}^h}(x) ~~; ~~
\mbox{ with } x\in \Gamma_{n}^h\}.
\end{equation}
Here $\vec{\eta}_{\Gamma^+_{n}}$ represents an approximation of the unitary normal vector of $\Gamma_n^h$  outward to $\Omega^+_n$. More specifically,  since $\Gamma_{n}^h$ is piecewise linear, its normal vector  $\vec{\eta}_{\Gamma_{n}^h}(x)$ is not well define when $x$ is a vertices of $\mathcal{T}^h_{n}$.  Different strategies can be used to handle this problem, for instance, we can define the normal vector as the average of the two adjacent normal vectors of $x$; or we can interpolate the vertices of $\Gamma_n^{h}$ by a smooth curve and define the normal vector of $\Gamma_n^h$ at $x$ as the normal vector of the smooth interpolation of $\Gamma_n^h$. In our numerics, we implemented the first strategy.

Next, the triangulation
$\mathcal{T}^h_{n+1}$ is defined such that $\Gamma^h_{n+1}$
is the union of edges in $\mathcal{T}^h_{n+1}$.
More precisely, we obtain $\mathcal{T}^h_{n+1}$ from
$\mathcal{T}^h_{n}$ using a harmonic
extension of the displacement $\tau\sigma^h \vec{\eta}_{\Gamma_n^h}$ as follows.
First, we introduce the vector function
$\vec{\omega}_n^h=(\omega_1^h,\omega_2^h)$ where
each component  satisfies
\begin{equation}\label{eq:current-triangulation-disp-plus}
\left\{ \begin{array}{ll}
\int_{\Omega^+_{n,h}} a_1\nabla \omega_j^h\nabla z^h\;dx=0& \mbox{ for all }
z^h\in {\cal P}^1(\mathcal{T}_n, \Omega^+_{n,h})\\
w^h_j(x^h)=0& \mbox{ for all }  x^h\in \Sigma^+ , ~\mbox{and}~x^h\in \mathcal{T}^h_n\\
w^h_j(x^h)=\tau \sigma^h(x^h)
(\vec{\eta}_{\Gamma^h_n}(x^h))\cdot
\vec{e}_j& \mbox{ for all  }   x^h\in \Gamma_n, ~\mbox{and}~x^h\in \cap \mathcal{T}^h_n,\\
\end{array}
\right.
\end{equation}
and
\begin{equation}\label{eq:current-triangulation-disp-minus}
\left\{ \begin{array}{ll}
\int_{\Omega^-_{n,h}} a_2\nabla \omega_j^h\nabla z^h\;dx=0& \mbox{ for all }
z^h\in {\cal P}^1(\mathcal{T}_n, \Omega^-_{n,h})\\
w^h_j(x^h)=0& \mbox{ for all }  x^h\in \Sigma^-, ~\mbox{and}~x^h\in \mathcal{T}^h_n\\
w^h_j(x^h)=\tau \sigma^h(x^h)
(\vec{\eta}_{\Gamma^h_n}(x^h))\cdot
\vec{e}_j& \mbox{ for all  }   x^h\in \Gamma_n, ~\mbox{and}~x^h\in \mathcal{T}^h_n,\\
\end{array}
\right.
\end{equation}
where $\vec{e}_1=(1,0)$ and $\vec{e}_2=(0,1)$. Then we define
the nodes of the new triangulation
\begin{equation}\label{eq:next-triangulation}
x_{n+1,j}^h=x_{n,j}^h+\vec{\omega}(x_{n,j}^h).
\end{equation}
The edges and triangles structures  of $\mathcal{T}^h_{n+1}$
is inherit directly from $\mathcal{T}^h_{n}$.

Finally, we observe that the initial  triangulation ${\cal T}_0$ can be  defined as any regular triangulation of $\Omega$  containing vertices  on the initial approximation $\Gamma_0$ of $\Gamma$.

\subsection{Summary of the iterative method}\label{fig:summary}
We now summarize the proposed iteration  for a given a tolerance
$\epsilon_{tol}$.
\begin{center}
\begin{minipage}{10cm}
Input: Domain $\Omega$ and boundary condition $g$.\\
Output: Free interface approximation, $\Gamma^h_n$,
and approximation of the  solution, $u^h_n$.
\begin{enumerate}
\item Set up $\Gamma_0$ (and the
positive and negative subsets $\Omega_0^+$ and
$\Omega_0^-$).
\item For $n=1,2,\dots$, until convergence, do
\begin{enumerate}
\item Compute $u_n^h$ by solving
(\ref{eq:current-u-plus-h})
and (\ref{eq:current-u-minus-h}).
\item Compute $\sigma_n^h$ in (\ref{eq:def:sigmah}).
\item Compute the triangulation displacement
$\vec{\omega}^h_n$ by solving (\ref{eq:current-triangulation-disp-plus})
and (\ref{eq:current-triangulation-disp-minus}).
\item Set up the new free interface approximation
$\Gamma_{n+1}^h$ in (\ref{eq:next-Gamma-h})
and triangulation $\mathcal{T}_{n+1}^h$ in (\ref{eq:next-triangulation}).
\end{enumerate}
\end{enumerate}
\end{minipage}
\end{center}
Here, the convergence criteria is given by $\|\sigma^h_n\|_{L^\infty(\Gamma)}< \epsilon_{tol}$.

\subsection{The parameter $\tau$}

In order to get some insight on the role of the parameter $\tau$ we may compare  our method with a regular   Newton's method to solve $\sigma(z^h)=0$, where $z^h$ represents the coordinates of the vertices of the partition belonging to the free boundary. Formally,   a  Newton's method for our problem would consists of the following the iteration
$$
\nabla_z^h \sigma( z^{n,h})  (z^{n+1,h}- z^{n,h}) =   -\sigma(z^{n,h}).
$$
Here $\nabla_z^h \sigma(z^{n,h})$ represents a formal  derivative operator of $\sigma$ with respect of $z^{n,h}$. Assuming it is possible to invert the operator
$\nabla_z^h \sigma(z^{n,h})$ we would have
$$
  (z^{n+1,h}- z^{n,h}) =   -(\nabla_z^h \sigma( z^{n,h}))^{-1} \sigma(z^{n,h})  .
$$
 From (\ref{eq:next-Gamma}) we conclude that our
method satisfies
$$
z^{n+1,h} - z^{n,h} =\tau \sigma(z^n) \vec{\eta}_{\Gamma^n}
$$
Finally,  assuming it is correct to update the $z^{n,h}$ by moving it toward the average of  normal directions of $\Gamma_n^h$ adjacent to $z^{n,h}$,  we expect to obtain   a damped Newton method by choosing $\tau$ sufficiently small.

In our numerics we observed that  the parameter $\tau$ should  be chosen  sufficiently small to avoid big variations of the  triangulation with respect to $\Omega$ in a single step.

\section{Applications  to jets of two fluids in a pipe}\label{sec:jet}

In this section we apply our method to a simplified version
of the jet problem for two fluids. We consider a  version of
 the model discussed in \cite{Alt-Caf-Fri_Jet-I}.
There, it is considered the model of  two planar flows
along an infinite pipe with  one free interface. Here we use our method to computed approximate solutions
in a truncated pipe  model with some given inflow/outflow data.

\begin{figure}[htb]
\centering
{\psfig{figure=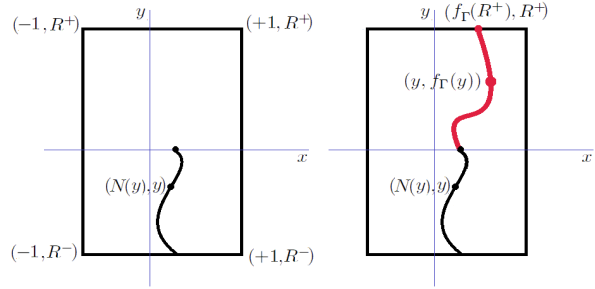,height=4cm,width=8cm,angle=0}}
\caption{Simple vertical  pipe configuration. }
\label{fig:jetgeom}
\end{figure}

The model  for jet flow studied in \cite{Alt-Caf-Fri_Jet-I} is the following. Let $u$ denote the stream function associated to the irrotational flow of two ideal fluids. The regions occupied by each different fluid are represented by  the support of $u^+$ and $u^-$, where  $u^+$ ($u^-$) denotes the positive (negative) part of $u$.  Let $N_1:\mathbb{R} \rightarrow (c_1,c_2)$, with $0<c_1<c_2$  be a continuous and piecewise $C^2$ function, satisfying
$\lim_{y\to \infty}N_1(y)=B$, $\int_0^\infty (N_1(y)-B)^2\;dy<\infty$, and $\int_0^\infty N_1'(y)^2\;dy<\infty$. The Alt {\it et. all.}  model assumes the two fluids occupy   an infinity semi-strip region enclosed by the graph  $\{(N_1(y),y); ~y>a, ~a<0\}$, and the lines $y=a$ and $x=-1$. The fluids enter the region at the boundary $\{y=a\}$, and the two fluids are separated  from each  other  in $\{y<0\}$ by a given continuous and piecewise $C^2$ curve $N_2: [a,0] \rightarrow (-1+\delta,c_2)$,  satisfying  $dist(N_2, N_1)>0$.  A special truncated  case of this configuration is shown in Figure \ref{fig:jetgeom} (left).  The problem consists in  finding the free boundary separating the  two fluids in the region $y>0$, assuming each flow has constant speed when $y \rightarrow \infty$.
More specifically, we look for
$u$ and $\lambda$ satisfying
\begin{equation}\label{eq:problema-orig-ACF}
 \begin{array}{l}
 \Delta u=0 \mbox{ in each fluid} \\
 |\nabla u^+|^2  - |\nabla u^-|^2= \lambda    \mbox{ on the free boundary separating the two fluids} \\
 u=Q  \mbox{ on  $\{x=-1\}$   and   $u=-1$ on  $N_1$   }\\
  u=g \mbox{ on $ [-1,N_1(a)]\times \{y=a\}$    }
\end{array}
\end{equation}
where
\begin{equation}\label{eq:def:lambda}
\left\{ \begin{array}{l}
\lambda=1/(1+b)^2 -Q^2/(B-b)^2, ~\mbox{and the free boundary} \\
\mbox{approximates the point $(b,0)$ when $y\to \infty$}.
\end{array}
\right.
\end{equation}
 Here the function $g \in C^1$ is monotone decreasing and
$$
 \begin{array}{lll}
0\leq g(x) \leq Q &\mbox{for}& x <N_2(a),\\
 -1\leq g(x) \leq 0 &\mbox{for}& x >N_2(a), \\
 ~g(-1)=Q  &\mbox{and}& g(N_1(a))=-1.
 \end{array}
 $$
Existence  and uniqueness of solution for this problem was studied in \cite{Alt-Caffarelli-Friedman}, where it was proved that minimizers of an appropriate functional are {\it weak} solutions of problem  (\ref{eq:problema-orig-ACF}).

We construct approximated solutions of the above problem. In
particular we work with a truncated domain to represent a pipe.

We now refer to the problem configuration
in Figure  \ref{fig:jetgeom}. Given a positive constant $Q$, and functions $N:[R^-,0]\to (-1,1)$
and $g:[-1,1]\to\mathbb{R}$,
we want to find $u:(-1,1)\times (R^-,R^+)\to \mathbb{R}$
and free interface $\Gamma$ represented by
\[
\Gamma=\{ (y,f_\Gamma(y)); \mbox{ with }  0\leq y\leq R^+ \}
\]
where $f_\Gamma:[0,R^+]\to \mathbb{R}$ is such that $f_\Gamma(0)=N(0)$  (see Figure \ref{fig:jetgeom} right picture).
The function $u$ and the free interface $\Gamma$ satisfy
\begin{equation}\label{eq:truncatedJet1}
\Delta u =0  \mbox{ in } \Omega^- \quad\mbox{ and } \quad \Delta u =0  \mbox{ in } \Omega^+
\end{equation}
where
\begin{eqnarray*}
&& \Omega^+=\{ -1<x<N(y),~ R^-<y\leq 0\}\cup \{-1<x<f_\Gamma(y),~ 0\leq y<R^-\}\\
&& \Omega^-=\{ N(y)<x<1,~ R^-<y\leq 0\}\cup \{f_\Gamma(y)<x<1,~ 0\leq y<R^-\}.
\end{eqnarray*}
The function $u$ has to satisfy the following known given data,
\begin{eqnarray}
&&u(N(y),y)=0,\quad  R^-<y<0;\\
&&u(x,R^-)=g(x),\quad  -1<x<1;\\
&&\frac{\partial u}{\partial y}u(x,R^+)=0,  \quad -1<x<1;\\
&&u(-1,y)=Q, \quad R^-<y<R^+;\\
&&u(1,y)=-1, \quad R^-<y<R^+;
\end{eqnarray}
and the following conditions on the free interface
\begin{eqnarray}
u=0, && \mbox{ on } \Gamma \quad (\mbox{ or } u(f_\Gamma(y),y)=0, 0\leq y<R^+\mbox{)}\\
(\partial_{\eta} u^{+} )^2-(\partial_{\eta} u^{-})^2=\lambda&& \mbox{ on } \Gamma
\end{eqnarray}
where $\lambda$ is given by
\begin{equation}\label{eq:truncatedJetFinal}
\lambda= \left(\frac{1}{1+b}\right)^2 -\left( \frac{Q}{2-b}\right)^2.
\end{equation}

We note that the boundary condition on the top of the domain (see Figure \ref{fig:jetgeom})
is the homogeneous Neumann boundary condition, hence, the free boundary
is not fixed at the top. That is, the value $f_{\Gamma}(R+)$ is not prescribed.

For each value of $b\in (-1,1)$, we can compute $\lambda$ through (\ref{eq:truncatedJetFinal}) and use our method
to find an approximation of  $u=u^b$  and $\Gamma=\Gamma^b$
(represented by $f_\Gamma=f_\Gamma^b$)  that
solve (\ref{eq:truncatedJet1})-(\ref{eq:truncatedJetFinal}).
Since the free boundary must  have the vertical line $x=b$ as an asymptote, a
feasible approximation of the free interface $\Gamma$  is obtained
if $b=b^*$ where
\begin{equation}\label{eq:condition-on-top}
b^*=f_\Gamma^{b^*}(R^+).
\end{equation}
 This is compatible
with the asymptote condition $\lim_{y\to \infty} f_\Gamma(y)=b^*$.

Next, we use  a {\it bisection} algorithm (applied to the function $(-1,1)\ni b\mapsto f_{\Gamma}^b(R^+)\in (-1,1)$)  to find the correct value
of $b^*$ such that (\ref{eq:condition-on-top}) is satisfied.

In the two examples presented next we run our method described in Subsection \ref{fig:summary} until
$\|\sigma^h_n\|_{L^\infty}<tol=10^{-6}$.

\begin{figure}[htb]
\centering
{\psfig{figure=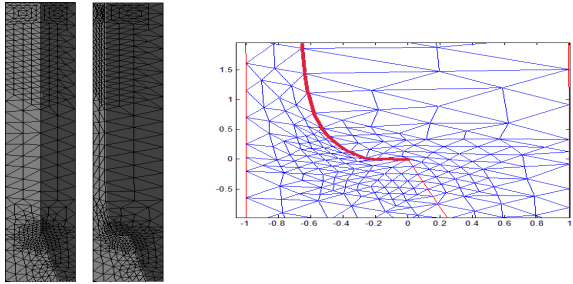,height=4cm,width=9cm,angle=0}}
\caption{Free interface that solves (\ref{eq:truncatedJet1})-(\ref{eq:truncatedJetFinal}) and
(\ref{eq:condition-on-top}) with
$N(y)=0.5|y|/R^-$, $g(x)=(x-0.5)/1.5$ if $-1<x<0.5$, $g(x)=5(x-0.5)/0.5$ if $0.5<x<1$,
 $Q=5$ and $b^*=1/3$.
Initial configuration  (Left), final configuration (Middle) and a zoom around
$(0,0)$ showing the free interface and the final mesh (Right).}
\label{fig:jet1}
\end{figure}
\begin{figure}[hbt]
\centering
{\psfig{figure=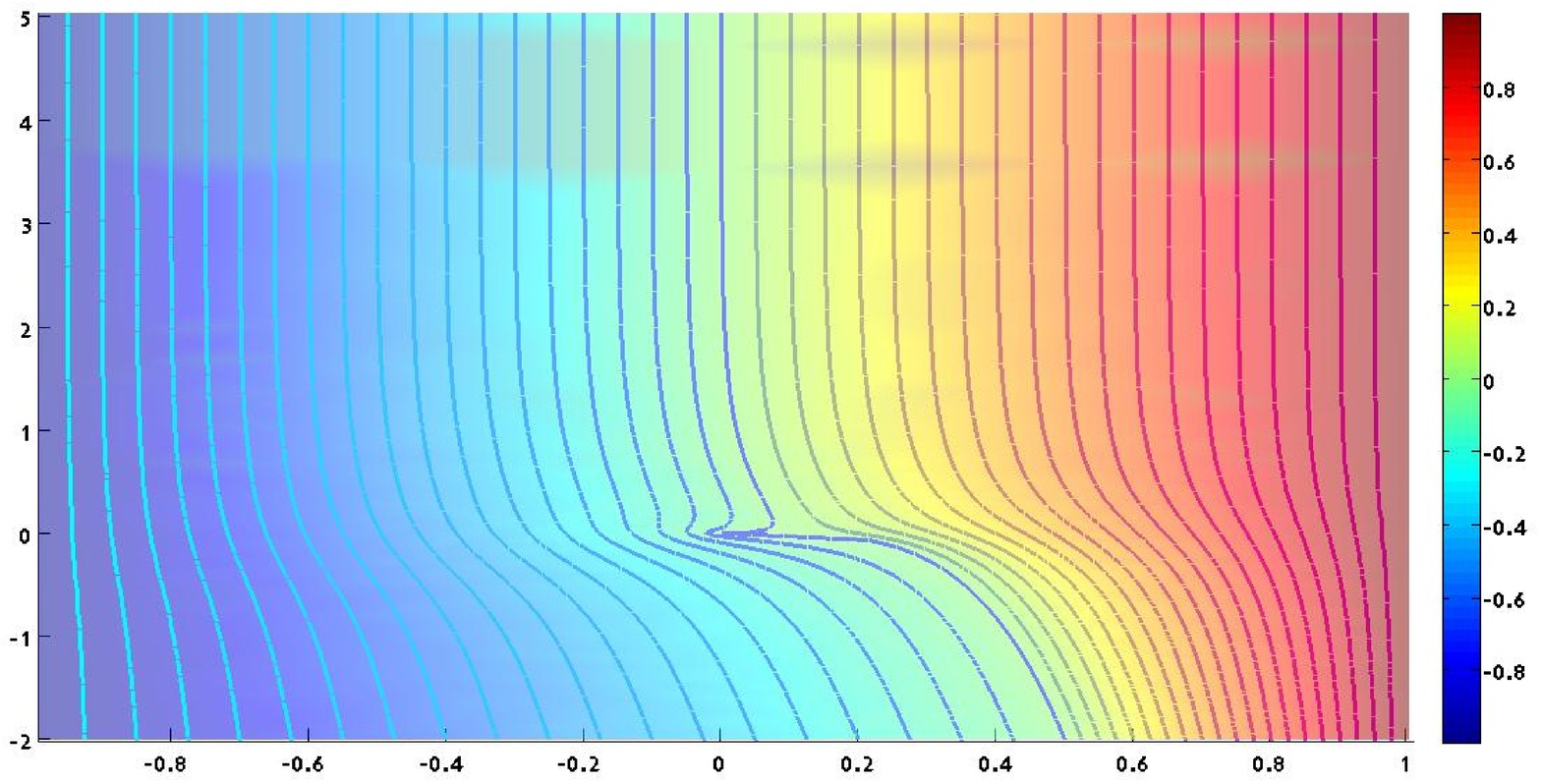,height=3cm,width=3cm,angle=0}}
{\psfig{figure=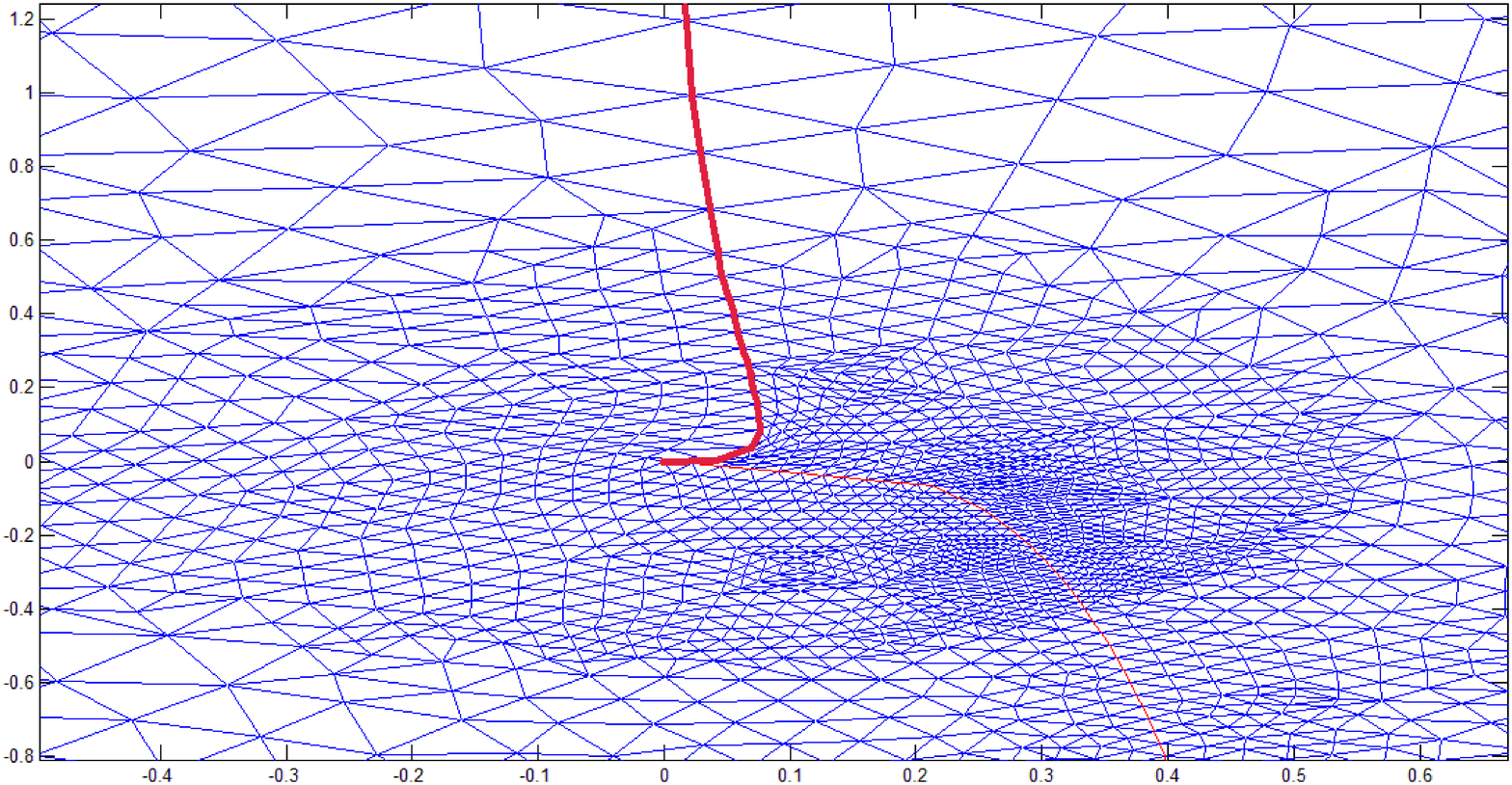,height=3cm,width=3cm,angle=0}}
\caption{ Function $u$ and free interface $\Gamma$ that
solve (\ref{eq:truncatedJet1})-(\ref{eq:truncatedJetFinal}) and (\ref{eq:condition-on-top}) with
 $N(y)=0.5(|y|/R^-)^{0.25}$,
$g(x)=(x-0.5)/1.5$ if $-1<x<0.5$ , $g(x)=(x-0.5)/0.5$ if $0.5<x<1$,
$Q=1$ and  $b^*=1$. Solution (Left) and  zoom around
$(0,0)$ showing the free interface and the final mesh (Right)}
\label{fig:jet2}
\end{figure}
The first example   considers  the nozzle represented by
$N(y)=0.5|y|/R^-$ with  $Q=5$.
The data on the bottom is given by,
$h(x)=(x-0.5)/1.5$ if $-1<x<0.5$
and $h(x)=5(x-0.5)/0.5$ if $0.5<x<1$.
We obtain $b^*=1/3$. The initial free boundary approximation $\Gamma_0$
is the strait line from $(0,0)$ to $(0,R^+)$
The resulting free boundary is displayed in Figure \ref{fig:jet1}.

In the second example of  jet flow problem we consider
the nozzle represented by $N(y)=0.5(|y|/R^-)^{0.25}$
with $Q=1$ and the Dirichlet data on the bottom side given by $h(x)=(x-0.5)/1.5$ if $-1<x<0.5$ and $h(x)=(x-0.5)/0.5$ if $0.5<x<1$. We obtained $b*=1$.
The initial free boundary approximation $\Gamma_0$
is the strait line from $(0,0)$ to $(0,R^+)$.
The resulting free boundary and numerical solution (with constant lines -stream lines)
are  displayed in Figure \ref{fig:jet2}.

\section{Application to a free boundary problem
arising in magnetohydrodynamics}\label{sec:mhd}

In this section we apply our methodology to the  model of \emph{plasma} problem
studied in \cite{friedman1, plasma2}.  Here we are interesting in modeling the plasma confined in a Tokamac machine. More specifically, given $\Omega\subset\mathbb{R}^2$ and the
positive constants $\gamma$ and $\lambda$, the plasma problem is to find $u$, a closed curve $\Gamma$
lying in $\Omega$ and a positive constant $\beta$ such that
\begin{equation}\label{eq:mhd}\left\{
\begin{array}{rll}
-\Delta u =\beta u && \mbox{ in } \Omega^-=\mbox{int}\{ x\in \Omega, u(x)\leq 0\},\\
\int_{\Omega^-} u^2 =1&& \\
-\Delta u=0 &&\mbox{ in } \Omega^+=\{ x\in \Omega, u(x)>0\},\\
u=\gamma &&\mbox{ on  } \partial \Omega,\\
u=0 &&\mbox{ on } \Gamma  \mbox{ and } \Gamma=\partial \Omega^-,\\
 |\partial_\eta u^+|^2  - |\partial_\eta u^-|^2= \lambda  &&\mbox{ on} ~\Gamma.
\end{array}\right.
\end{equation}
Here the plasma is enclosed by the curve $\Gamma$, and the complement of this region with respect to $\Omega$ is vacuum.  The function $u$ represents a flux function associated to the magnetic induction $\vec{B}$,  satisfying
$ \vec{B}=(u_{x_2},-u_{x_1},0)$.

It is easy to modify our method and apply it to this problem.
We follow the description in Subsection \ref{fig:summary} and iterate until
$\|\sigma^h_n\|_{L^\infty}<tol=10^{-6}$.
In this  problem, the free boundary is a closed curve separating the domain in two connected components; as showed in \cite{friedman1}. The adaptation of our scheme to treat this problem is straightforward.
We also mention that other formulations of the model, having
the nonlinear condition on the free boundary, are also possible; see \cite{friedman1} and references therein.

\begin{figure}[htb]
\centering
{\psfig{figure=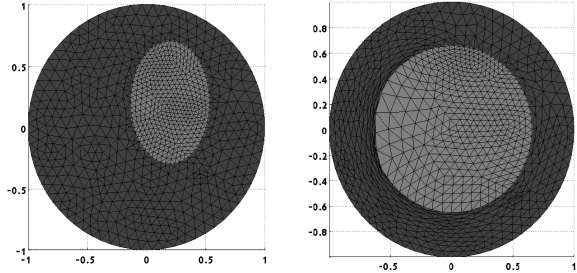,height=4cm,width=8cm,angle=0}}
\caption{Free interface that solves (\ref{eq:mhd}) with
$\Omega=\{(x,y) : x^2+y^2<1\}$ , $\gamma=1$ and $\lambda=2^2-1^2=3$.
Initial configuration: an ellipse
centered at $(1/5,1/5)$ and with  axis $1/3$ and $1/2$. (left). Final configuration
(right).}
\label{fig:plasma1}
\end{figure}

In the first example, we consider
the case of $\Omega$ being the ball with center $(0,0)$ and radius 1.
We choose $\gamma=1$ and $\lambda=2^2-1^2=3$. The initial approximation
of the free boundary is  an ellipse
centered at $(1/5,1/5)$ and with  axis $1/3$ and $1/2$.
The resulting configuration is depicted in Figure \ref{fig:plasma1} and
$\beta=13.6727$.
We observe that the final shape of the free boundary approximates a circular region.
This coincide with the results in \cite{plasma2} where the authors proved
that for the domain $\Omega$ being the unit  circle, the resulting free boundary is circular and centered at $(0,0)$.
We note that in this example, the initial approximation of the free interface is far-off from the
solution and despite of this fact, our algorithm still converges to the solution.

\begin{figure}[htb]
\centering
{\psfig{figure=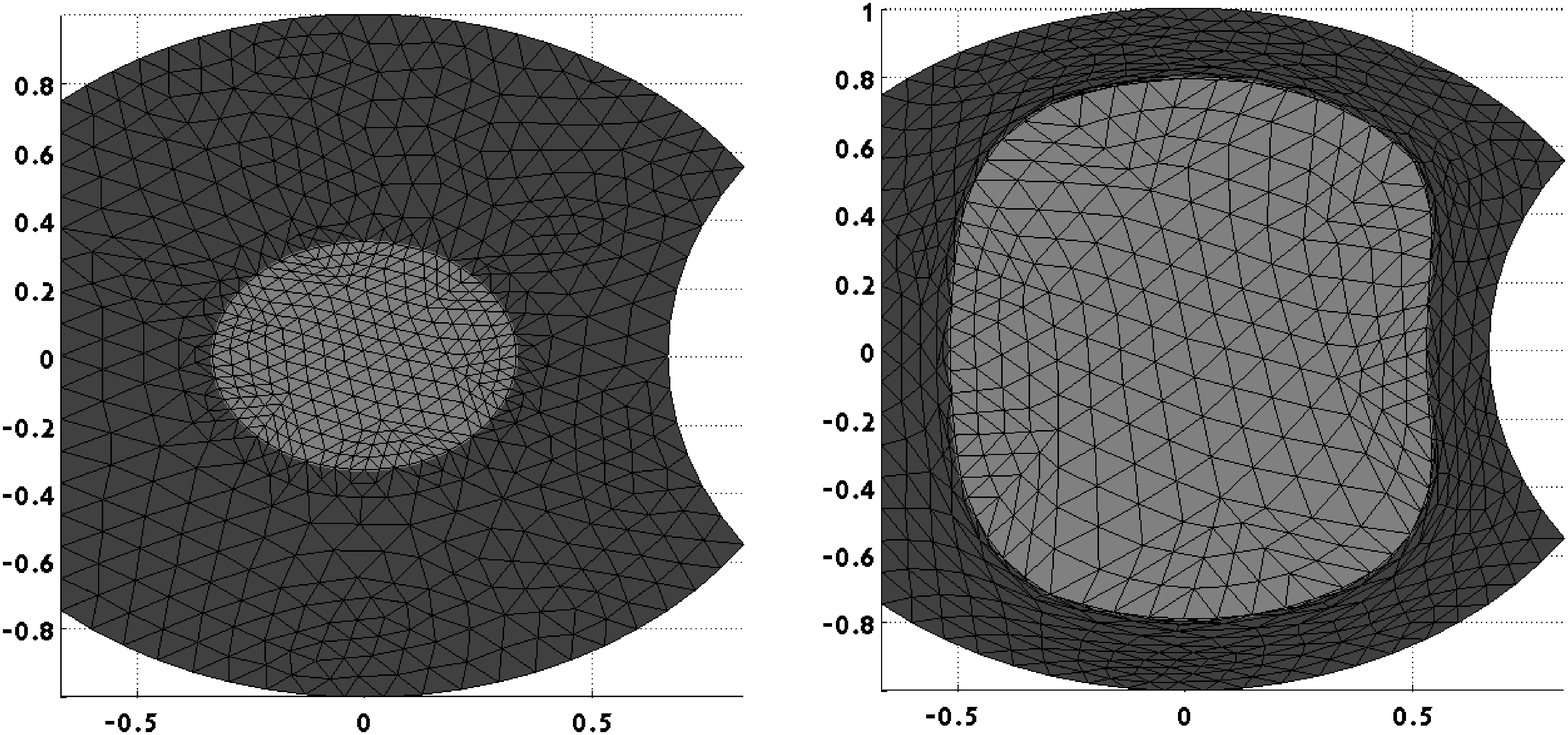,height=3cm,width=6cm,angle=0}}
{\psfig{figure=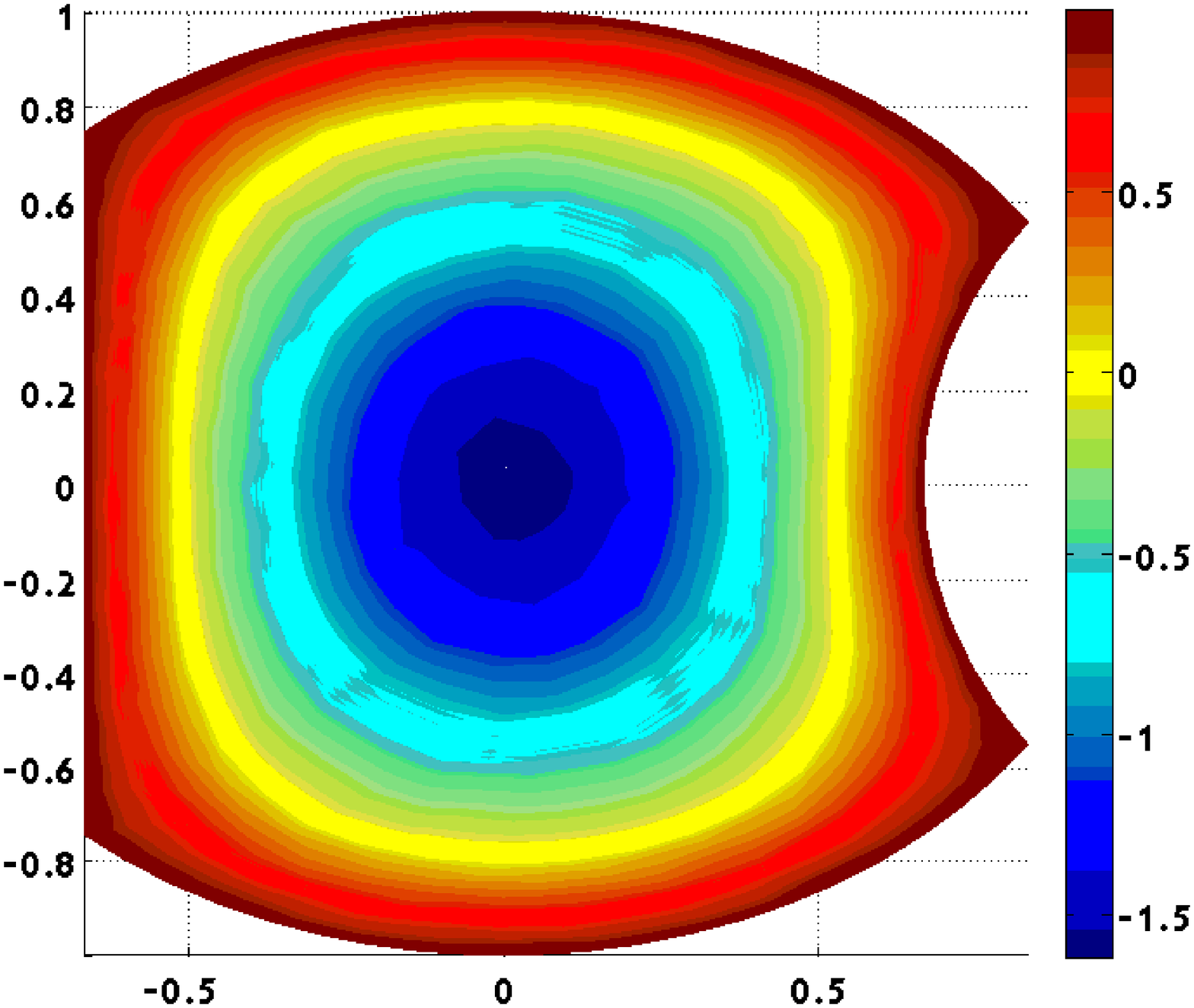,height=3cm,width=3cm,angle=0}}
\caption{Free interface that solves (\ref{eq:mhd}) with depicted
domain $\Omega$, $\gamma=1$ and $\lambda=5^2-1^2=4$.
Initial configuration: the circle with center $(0,0)$ and radius $1/3$ (Left). Final configuration
(Center). Solution (Right)}
\label{fig:plasma2}
\end{figure}

The second  example  considers the configuration described  in Figure \ref{fig:plasma2}.
The domain $\Omega$ corresponds to a circle from which it have been cut-off
the region $\{y<-2/3\}$ and the intersection with the circle with center $(5/3,0)$ and radius 1.
In this example we use $\lambda=5^2-1^2=4$.
The initial approximation of the free boundary is a ball with center $(0,0)$
and radius $1/3$.
The resulting free boundary is presented in Figure \ref{fig:plasma2} (center) and
the solution is plotted in Figure \ref{fig:plasma2} (right).
The computed value of $\beta=13.7034$.

\section{Additional numerical examples}\label{sec:additional-num}

In this section we present some representative numerical examples.
We run our method described in Subsection \ref{fig:summary} until
$\|\sigma^h_n\|_{L^\infty}<tol=10^{-6}$.

\subsection{A known free interface and error decay}\label{sec:known}
We consider problem  (\ref{eq:differential-problem}) with $a_1=a_2=1$ and  a known  exact solution, what allows us to measure the accuracy of our method.
The domain is
$\Omega=[0,1]\times [0,1]$,
 $\lambda=-1$, the boundary data is given by $g(x,y)=2\min\{x-0.5,0\}+\max\{x-0.5,0\}$.
 We note that $u(x,y)=2\min\{x-0.5,0\}+\max\{x-0.5,0\}$
is a solution of the problem. For this exact solution, the  free interface
is the strait vertical line from $(0.5,0)$ to $(0.5,1)$.

\begin{figure}[htb]
\centering
{\psfig{figure=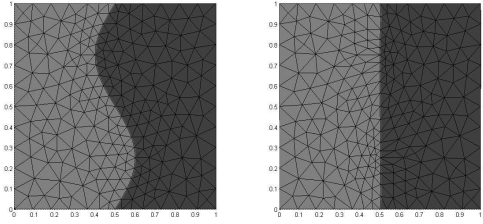,height=4cm,width=8cm,angle=0}}
\caption{Results for the test problem in
Section \ref{sec:known} Initial configuration (left). Final configuration
(right). }
\label{fig:line}
\end{figure}

\begin{figure}[htb]
\centering
{\psfig{figure=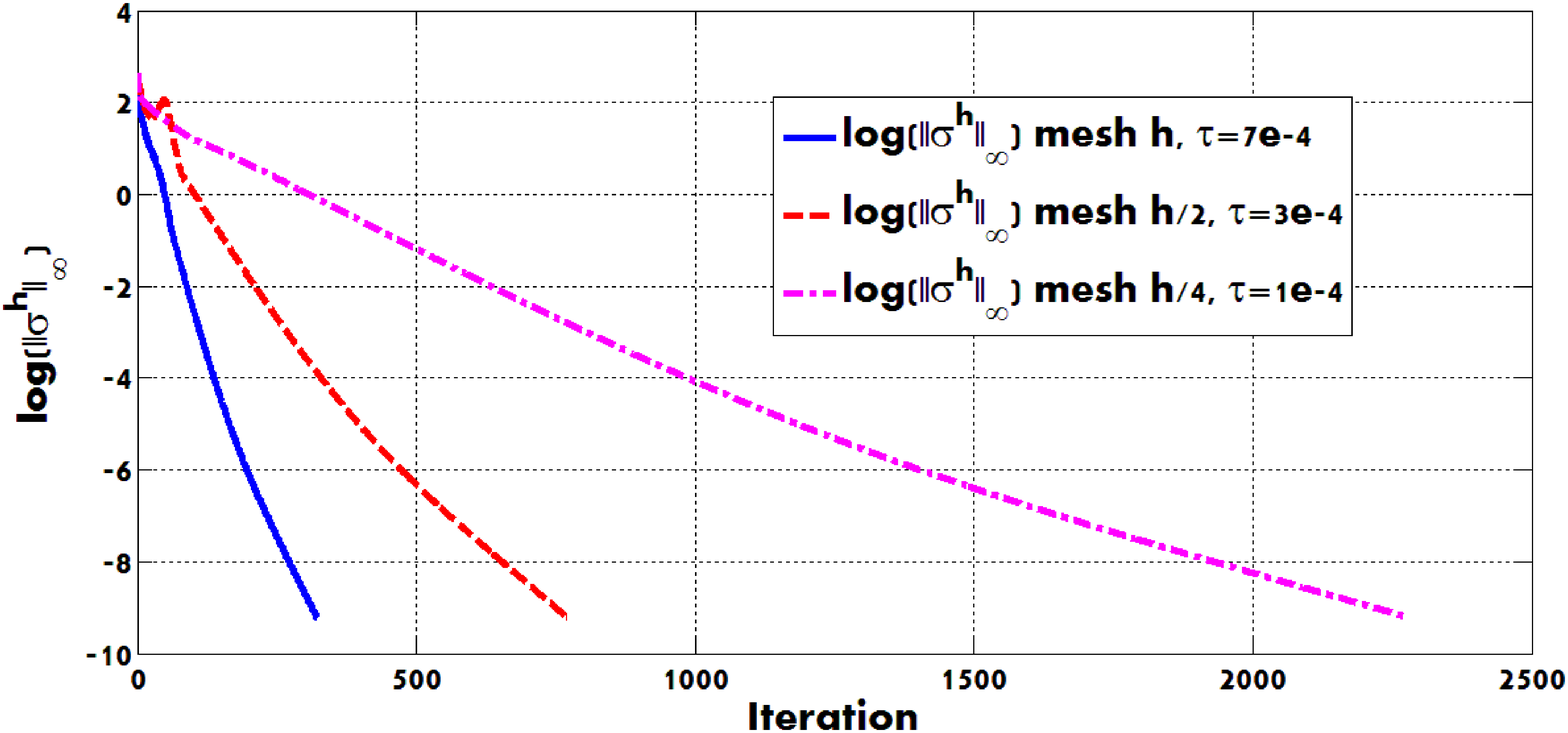,height=3.3cm,width=6cm,angle=0}}
{\psfig{figure=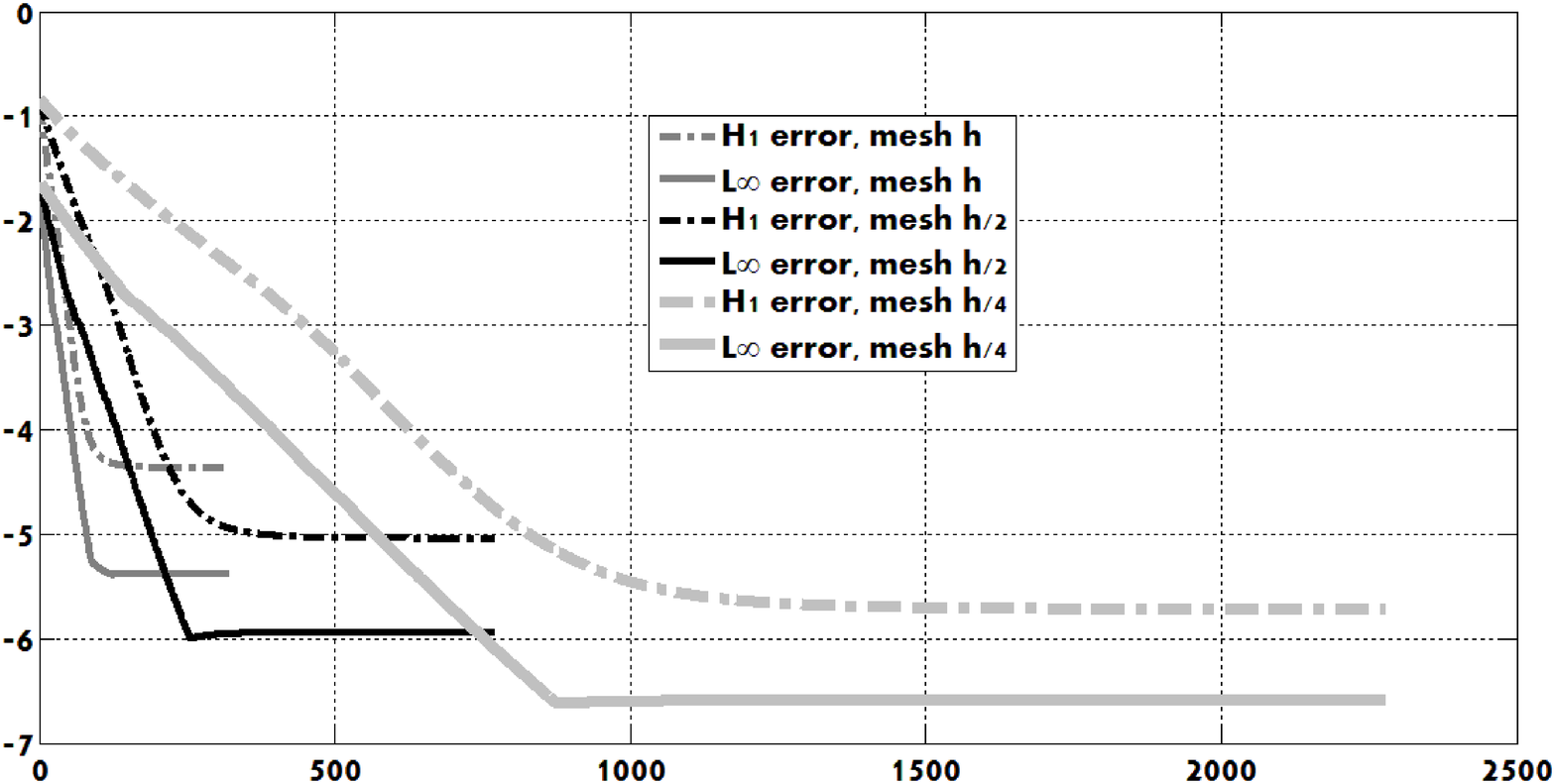,height=3cm,width=6cm,angle=0}}
\caption{ Result for the test problem in Section \ref{sec:known}.
We plot $\log(\|\sigma^h\|_{\infty})$ over the iterations
requires to achieve $\|\sigma^h\|_{\infty}<10^{-4}$ for
 three different
meshes (left). Corresponding solution error (in log scale) in
 $| \cdot|_{H^1(D)}$ and $L^\infty(D)$ norms.  }
\label{fig:sigmavalue1}
\end{figure}

We apply our method with the  initial approximation of the free boundary
 given by $\Gamma_0=\{(x,y) ; y=0.5+0.1\sin(2\pi y)\}$
and the parameter $\tau=10^{-4}$.
We present the initial and final subdomain
configuration in Figure \ref{fig:line}. In Figure
\ref{fig:sigmavalue1} we present the $L^\infty$ norm
of $\log(|\sigma^h|)$ in (\ref{eq:def:sigmah}) along the number of iterations $n$.
We observe a decay of the value $|\sigma^h|$ faster than $O(e^{-0.004n})$.
This example also show that our stopping criteria is effective in the sense
that at the last iteration, we see that the $L^\infty$ norm is already
in stagnated plateau for the corresponding mesh size.

\subsection{An example with heterogeneous coefficients}\label{sec:heterogeneous}
\begin{figure}[htb]
\centering
{\psfig{figure=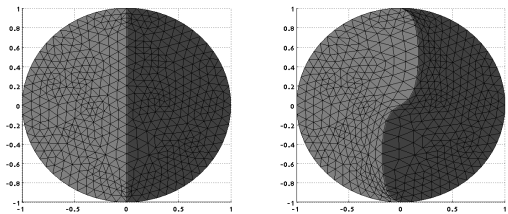,height=3cm,width=7.5cm,angle=0}}
\caption{Results for the test problem in
Section \ref{sec:heterogeneous} that involves heterogeneous coefficients. Initial configuration (left). Final configuration
(right) }
\label{fig:circle}
\end{figure}

This example considers a problem of type (\ref{eq:differential-problem})  with
heterogeneous coefficient in each side of the free interface.
We consider $\Omega=\{(x,y) : x^2+y^2<1\}$ and the coefficient(s)
\[
a_1(x,y)=\left\{\begin{array}{cc}
100,&  y\geq 0\\
1, & y<0
\end{array}\right.\quad \mbox{ and } \quad
a_2(x,y)=\left\{\begin{array}{cc}
1,&  y\geq 0\\
100, & y<0.
\end{array}\right.
\]
The Dirichlet data around the circle is given by  $g(x,y)=x$ and we use $\lambda=-1$.
The initial approximation of the free boundary is the strait line $\Gamma_0=\{ (0,y), 0\leq y\leq 1\}$.
We run our method with $\tau=10^{-5}$. We show the resulting free boundary
 in Figure \ref{fig:circle}.

\section{Conclusions and  comments}\label{sec:conclusions}
We have proposed a simple iterative method to handle free boundary problems
involving nonlinear  flux conditions. It is important to note, that the numerical treatment
of nonlinear flux conditions on the free interface have not been extensively
studied in the literature. This is the case despite of the fact that
the mathematical analysis of simple models with nonlinear
flux conditions on the free interface have been carried out by
Caffarelli and coauthors a couple  of decades ago.
The proposed  method is a simple  domain decomposition method
with inexpensive iterations. As a consequence it can be used for the better understanding
of simplified models of complex flow problems.
We present numerical results  showing that, our iterative method is effective
and perform well in several applications where
nonlinear flux jump constrain drive the free interface behavior.

We obtained encouraging numerical results with our method  but,
 its  mathematical analysis  is still needed.
In a future work we plan to address mathematically questions
related to the converge of the method.
Other  interesting numerical  aspects
we want to address are related to the implementation of adaptive refinement, the
use of inexact local solvers (instead of exact subdomain solvers),
and  the design  of  preconditioners for our scheme.
The extension to three dimensions can be considered.

We note that we consider simplified models of complicated flow problems. If
we want to extend our method for more realistic models we  need  to consider
  time dependent problems. In this case, it would be important to be able to handle topological changes in the evolution of the free boundary.

\section*{Acknoledgements}
The authors are thankful to Prof. Eduardo Teixeira for bringing this problem to our attention.
H.M.V. was partially supported by FAPERJ grants E-26/102.965/2011,
E-26/111.416 /2010.
J.G. research is based in part on work supported by Award No. KUS-C1-016-04, made by King Abdullah University of Science and Technology (KAUST).
\appendix
\section{An approximation of the flux}\label{sec:flux}

Given a free interface approximation
$\Gamma^h$, we consider the approximation
of the flux of $u^h$ (the solution  of   problem (\ref{eq:current-u-plus-h})) on ${\Omega^+_{h}}$.

Denote by $A=[a_{ij}]$ the Neumann finite element matrix
defined by
\[
a_{ij}=\int_{\Omega^+_h}a_1\nabla \phi_{i}\nabla\phi_{j}\;dx
\]
where $\{\phi_i\}$ are the usual hat basis function
of the space ${\cal P}^1(\mathcal{T}^h_n, \Omega^+_{h})$.

We classify the nodes in interior nodes I,
boundary nodes $\Sigma^+$ and interface notes $\Gamma$.
This classification gives the following block structure
of the matrix $A$,
\[
A=
\left(\begin{array}{ccc}
A_{II} & A_{I\Sigma^+} & A_{I\Gamma}\\
A_{I\Sigma^+}^T & A_{\Sigma^+\Sigma^+} & A_{\Sigma^+\Gamma}\\
A_{I\Gamma}^T & A_{\Sigma^+\Gamma}^T & A_{\Gamma\Gamma}\\
\end{array}
\right).
\]

The solution of (\ref{eq:current-u-plus-h}) is given by,
\[
u^h=
\left(\begin{array}{c}
u_{I}^h\\
u_{\Sigma^+}^h\\
u_{\Gamma}^h\\
\end{array}
\right)=
\left(\begin{array}{c}
A_{II}^{-1}(A_{I\Sigma^+}g^h)\\
g^h\\
0\\
\end{array}
\right).
\]

We define $\mu$ by
\[
\mu=A_{I\Gamma}^T u_{I}=A_{I\Gamma}^TA_{II}^{-1}A_{I\Sigma^+}g^h.
\]

Let $N_\Gamma$ be the number of vertices of $\cal{T}^h$ on $\Gamma^h$.  We note that, using basic finite element analysis, we
see that $\mu=(\mu_{i}) \in \mathbb{R}^{N_\Gamma}$   with
\[
\mu_{i} =  \int_{\Omega^+_{h}}  (a_1 \nabla u^h ) \cdot \nabla \phi_{\ell_i} \;dx  = \int_{\Gamma^h}  (a_1 \nabla u^h )\cdot \eta_{\Gamma^h} \phi_{\ell_i} \;ds.
\]
Here given  $i\in \{1,...N_\Gamma\}$,   $\ell_i$ represents the index of the a node of $\mathcal{T}^h$ belonging to  $\Gamma^h$.

We use $\mu$  to obtain  a piecewise linear approximation of the flux
$ \nabla u^h \cdot \eta_{\Gamma^h}$. Since $u^h=0$ on $\Gamma^h$, for each edge of $e_k$ of $\Gamma^h$ we have
$$
\nabla u^h|_{e_k}= \partial_{\eta_k} u   \eta_k
$$
where $\eta_k$ represents the normal vector to edge  $e_k$ pointing in the outward direction of $\Omega^+_{h}$. Hence,
\begin{equation}\label{eq:mu_i}
\mu_{i} =  \int_{\Gamma^h}   (\eta_{\Gamma^h}^T  a_1  \eta_{\Gamma^h}) \partial_{\eta_{\Gamma^h}} u  \phi_{\ell_i} \;ds.
\end{equation}

We define $\lambda^h_1$ the  piecewise linear approximation  of  $\partial_{\eta_{\Gamma^h}} u$  as follows. First, we observe that  $\lambda^h_1 \in span \{ \phi_{\ell_i}|_{\Gamma^h}\}_{1 \leq i \leq N_\Gamma }$. Next   we introduce  the matrix
$Q=[q_{ij}]   \in \mathbb{R}^{N_\Gamma \times N_\Gamma}$ with
\[
 q_{i j}=\int_{\Gamma^h}  (\eta_{\Gamma^h}^T  a_1  \eta_{\Gamma^h}) \phi_{\ell_i}\phi_{\ell_j} \;ds.
\]
Finally, based on relation (\ref{eq:mu_i})  we define
\begin{equation}\label{eq:def_coef_lambda}
\lambda^h_1= \sum_{i}^{N_\Gamma} \alpha_i \phi_{\ell_i}|_{\Gamma^h}
\end{equation}
where $\alpha=(\alpha_i)$ is the solution of
\[
Q \alpha=\mu.
\]

In a similar way we define
$\lambda_2^h$, the approximation of
of the flux on $\Gamma^h$, of the solution of
(\ref{eq:current-u-minus-h}).

\begin{remark}
A more regular approximation of the flux can be done in practice. For instance, we could obtain $\alpha$ as the solution of the following problem
\[
(Q +\epsilon D) \alpha=\mu.
\]
where $D$ is diffusion of operator on $\Gamma^h$ and $\epsilon$ is a regularization parameter.
\end{remark}

\bibliographystyle{plain}
\bibliography{references}



\end{document}